\documentclass[12 pt]{report}
\usepackage{amssymb,latexsym,epsfig,mathrsfs,graphpap,graphics,amsmath,amssymb}
\usepackage{amsmath}

\begin{document}

\vspace{.2in}\parindent=0mm

\begin{flushleft}

     {\bf\Large {Fractional Multiresolution Analysis  and   \vspace{.1in} Associated Scaling Functions   in  $L^2(\mathbb R)$}}

  \parindent=0mm \vspace{.3in}
{\bf{ Owais Ahmad$^{*}$, Neyaz A. Sheikh$^{1}$  and  Firdous A. Shah$^{2}$ }}

\end{flushleft}

\parindent=0mm \vspace{.1in}
{{\it $^{*}$Department of  Mathematics,  National Institute of Technology, Hazratbal, Srinagar -190 006, Jammu and Kashmir, India. E-mail: $\text{siawoahmad@gmail.com}$}}

\parindent=0mm \vspace{.1in}
{{\it $^{1}$Department of  Mathematics,  National Institute of Technology, Hazratbal, Srinagar -190 006, Jammu and Kashmir, India. E-mail: $\text{neyaznit@yahoo.co.in}$}}

\parindent=0mm \vspace{.1in}
{{\it\small$^{2}${Department of  Mathematics,  University of Kashmir, South Campus, Anantnag-192 101, Jammu and Kashmir, India. E-mail: $\text{fashah79@gmail.com}$}}

\parindent=0mm \vspace{.2in}
{\bf{Abstract:}}  In this paper, we show  how to construct an orthonormal basis from Riesz basis by assuming that the fractional translates of a single function in the core subspace of the fractional multiresolution analysis form a Riesz basis instead of an orthonormal basis. In the definition of fractional multiresolution analysis, we show that the intersection triviality condition follows from the other conditions. Furthermore, we show that the union density condition also follows under the assumption that the fractional Fourier transform of the scaling function is continuous at $0$. At the culmination, we provide the complete characterization of the scaling functions associated with fractional multiresolutrion analysis.
 
\parindent=0mm \vspace{.2in}
{\bf{Keywords:}}}   Fractional MRA; Fractional Fourier transform;  Scaling function. 

\vspace{.1in}\parindent=0mm


\parindent=0mm \vspace{.2in}
{\bf{1. Introduction}}

\parindent=0mm \vspace{.1in}
Fourier transform  is one of the most valuable and frequently used tools in signal processing and analysis. For Fourier transform, a signal can be represented either in the time or in the frequency domain, and it can be viewed as the time-frequency representation of a signal. In 1980, Victor Namias \cite{10} introduced the concept of fractional Fourier transform (FrFT) as a generalization of the conventional Fourier transform to  solve certain problems arising in quantum mechanics.  It is also referred as {\it rotational Fourier transform} or {\it angular Fourier transform} since it depends on a parameter $\alpha$ which is interpreted as a rotation by an angle $\alpha$ in the time-frequency plane. Like the ordinary Fourier transform  corresponds to a rotation in the time frequency plane over an angle $\alpha = 1\times {\pi}/{2}$, the FrFT corresponds to a rotation over an arbitrary angle $\alpha= \rho\times {\pi}/{2}$ with $\rho\in\mathbb R$.It has applications in different fields like quantum mechanics \cite{10}, optics \cite{{a1},{a2}}, signal processing \cite{{b1},{b2},{b3},{b4},{b5}}, and image processing \cite{{c1},{c2},{c3}}. Although the FrFT has a number of attractive properties, the fractional Fourier representation of a signal only provides overall FrFD- frequency content with no indication about the occurrence of the FrFD spectral component at a particular time. Since the FrFT uses a global kernel like Fourier transform, it fails in locating the FrFD spectral contents which is required in some applications. The concept of FrWT was initially proposed in \cite{men}, where FrFT is firstly used to derive the fractional spectrum of a signal and wavelet transform is then performed on the obtained fractional spectrum. Since the fractional spectrum derived by the FrFT only represents the FrFD-frequency over the entire duration of the signal, the FrWT defined in \cite{men} actually fails in obtaining the information of the local property of the signal. In \cite{3}, a fractional wave packet transform was developed and the basic idea is to introduce the wavelet basis function to FrFT. More recently, a new FrWT was proposed in \cite{shi} based on the concept of fractional convolution.  In \cite{ofwpf},  the notion of fractional wavepacket systems in $L^2(\mathbb R)$ is introduced and the correponding frames are characterized. 

\parindent=8mm \vspace{.1in}
Multiresolution analysis is an important mathematical tool since it provides a natural framework for understanding and constructing discrete wavelet systems.  The concept of MRA has been extended in various ways in recent years. These concepts are generalized to  $L^2\big(\mathbb R^d\big)$, to lattices different from  $\mathbb Z^d$, allowing the subspaces of MRA to be generated by Riesz basis instead of orthonormal basis, admitting a finite number of scaling functions, replacing the dilation factor 2 by an integer $M\geq 2$ or by an expansive matrix $A\in GL_{d}(\mathbb R)$ as long as $A\subset A\mathbb Z^d$. All these concepts are developed on regular lattices, that is the translation set is always a group. In the heart of any MRA, there lies the concept of  scaling functions. Cifuentes et al.\cite{ref3} characterized the scaling function of MRA in a general settings .The multiresoltion analysis  whose scaling functions are characteristic functions  some elementary properties of MRA of $L^2(\mathbb{R}^n)$  are established by Madych \cite{ref11}. Zhang  \cite{ref17} studied  scaling functions of standard MRA and wavelets. Zhang \cite{ref17} characterized support of the Fourier transform of scaling functions. Malhotra and Vashisht \cite{lalit} provides the characterization of scaling functions on Euclidean spaces. The multiresolution analysis (MRA) associated with corresponding to  FrWT \cite{shi} was then given in \cite{2}. Since this kind of FrWT analyze the signal in time-frequency-FrFD domain, its physical meaning requires deeper interpretation. Another kind of FrWT which was developed in \cite{ap} solves the issue in \cite{shi} since the analysis only involves time-FrFD domain. However, the MRA associated with this kind of FrWT is not addressed. The main objectives of this article are as follows:
\begin{itemize}
\item To show how to construct an orthonormal basis from Riesz basis by assuming that the fractional translates of a single function in the core subspace of the fractional multiresolution analysis form a Riesz basis instead of an orthonormal basis.
\item To derive intersection triviality condition in the definition of fractional multiresolution analysis from other conditions.
\item To prove that union density condition  follows under the assumption that the fractional Fourier transform (FrFT) of the scaling function is continuous at $0$.
\item To provide the characterization of the scaling functions associated with fractional multiresolution analysis.
\end{itemize}

\parindent=8mm \vspace{.0in}
The rest of the article is structured as follows. In section 2, we discuss preliminaries of fractional Fourier and fractional wavelet transforms. Section 3 is devoted to the construction of an orthonormal basis from Riesz basis. In section 4, we show that the properties in the defintion of fractional multiresolution analysis are not independent. We show that the intersection triviality condition in the definition of fractional multiresolution analysis from other conditions and union density condition  follows under the assumption that the fractional Fourier transform (FrFT) of the scaling function is continuous at $0$. In section 5, we characterize the scaling functions associated with fractional multiresolution analysis.

\parindent=0mm \vspace{.2in}
{\bf{2. Preliminaries}}

\parindent=0mm \vspace{.1in}
This section gives the basic background to the theory of fractional Fourier and wavelet transforms which is as follows.

\parindent=0mm\vspace{.1in}
The fractional Fourier transform with parameter $\alpha$  of function $f(t)$ is defined by

$$\mathcal{F}_\alpha\big\{f(t)\big\}(\xi)=\hat f^\alpha(\xi)=\int_{-\infty}^{\infty} {\mathcal K}_\alpha (t,\xi)f(t)\,dt,\eqno(2.1)$$

\parindent=0mm \vspace{.1in}
where $  {\cal K}_\alpha (t,\xi)$ is called kernel of the FrFT given by
$${\cal K}_\alpha (t,\xi)=\left\{\begin{array}{ll}
  C_\alpha \exp\Big\{i(t^2+\xi^2)\dfrac{\cot\alpha}{2}-it\xi \csc\,\alpha\Big\},&\alpha\neq n\pi, \\
\delta(t-\xi), &\alpha=2n\pi, \\
\delta(t+\xi), &\alpha=(2n\pm 1)\pi, \\
\end{array}{}\right.\eqno(2.2)$$

\parindent=0mm \vspace{.1in}
$\alpha={\rho\pi}/{2}$ denotes the rotation angle of the transformed signal for FrFT,   the FrFT operator is designated by $\mathcal{F}_\alpha$ and
\begin{align*}
C_\alpha=\left(2\pi i \sin\alpha\right)^{-1/2}e^{i\alpha/2}=\sqrt{\dfrac{1-i\cot\alpha}{2\pi}}.\tag{2.3}
\end{align*}

\parindent=0mm\vspace{.0in}
The corresponding inversion formula is given by
$$f(t)=\int_{-\infty}^{\infty}\overline{{\cal K}_{\alpha}(t,\xi)}\,\hat{f}^\alpha(\xi)\,d\xi,\eqno(2.4)$$
where
\begin{align*}
  \nonumber {\cal K}_{\alpha}(t,\xi) &= \frac{(2\pi i\sin\alpha)^{1/2}\,e^{-i\alpha/2}}{\sin\alpha}\cdot \exp\left\{\frac{-i(t^2+\xi^2)\cot\alpha}{2}+it\xi \csc\,\alpha\right\}\\
  \nonumber &=\overline{C_\alpha}\exp\left\{\frac{-i(t^2+\xi^2)\cot\alpha}{2}+i t\xi \csc\,\alpha\right\}\\
   &={\cal K}_{-\alpha}(t,\xi)\tag{2.5}
\end{align*}
and
$$C_{\alpha}=\frac{(2\pi i\sin\alpha)^{1/2}e^{-i\alpha/2}}{2\pi \sin\alpha}=\sqrt{\dfrac{1+i\cot\alpha}{2\pi}}=C_{-\alpha}.\eqno(2.6)$$

\parindent=0mm \vspace{.0in}
{\bf {Definition 2.1.}} A  fractional wavelet is a function $\psi\in L^2(\mathbb{R})$ which satisfies the following condition:
$$C_{\psi}^{\alpha}=\int_{\mathbb{R}}\frac{\left|\mathcal{F}_\alpha\left\{e^{{-i(t-\xi)^2}/{2}\,\cot\alpha}\psi\right\}(\xi) \right|^2}{|\xi|}\,d\xi<\infty,\eqno(2.7)$$

\parindent=0mm\vspace{.1in}
where $\mathcal{F}_\alpha$ denotes the FrFT operator.

\parindent=8mm\vspace{.1in}
Analogous to the classical wavelets, the fractional wavelets can be obtained from a fractional mother wavelet $\psi\in L^2(\mathbb{R})$ by the combined action of translation and dilations as
$$\psi_{{\alpha},a,b}(t)=\frac{1}{\sqrt{a}}\,\psi\left(\frac{t-b}{a}\right)\exp\left\{\frac{-i(t^2-b^2)\cot\alpha}{2}\right\}\eqno(2.8)$$

\parindent=0mm\vspace{.1in}
where $a\in\mathbb R^+$ and $b\in\mathbb R$ are scaling and translation parameters, respectively. 
If $\alpha={\pi}/{2}$, then $\psi_{\alpha, a,b}$ reduces to the conventional wavelet basis.

\parindent=5mm\vspace{.1in}
Note that if $\psi(t)\in L^2(\mathbb{R})$, then $\psi_{\alpha, a,b}(t)\in L^2(\mathbb{R})$,
$$\left\|\psi_{\alpha, a,b}\right\|^2_{2}=|a|^{-1}\int_{-\infty}^{\infty}\left|\psi\left( \frac{t-b}{a}\right)\right|^2dt = \int_{-\infty}^{\infty}\big|\psi(y)\big|^2dy =\big\|\psi\big\|_{2}^{2}.$$

\parindent=0mm\vspace{.1in}
Moreover, the fractional Fourier transform of $\psi_{\alpha, a,b}(t)$ is given by

\parindent=0mm\vspace{.1in}
$\mathcal{F}_\alpha\big\{\psi_{\alpha, a,b}(t)\big\}$
$$=\sqrt{a}\,\exp\left\{\frac{i(b^2+\xi^2)\cot\alpha}{2}-ib\,\xi \csc\,\alpha-\frac{ia^2\xi^2\cot\alpha}{2}\right\}\mathcal{F}_\alpha\Big\{ e^{{-i(\cdot)^2\cot\alpha}/{2}}\psi\Big\}(a\xi)\eqno(2.9)$$

\parindent=5mm\vspace{.1in}
The continuous fractional wavelet transform (FrWT) of  function $f\in L^2(\mathbb R)$ with respect to an analyzing wavelet $\psi\in L^2(\mathbb R)$ is defined as
$${\mathscr W}^{\alpha}_{\psi}f(a,b)=\big\langle f,\psi_{\alpha, a,b}\big\rangle= \frac{1}{\sqrt{a}}\int_{-\infty}^{\infty}f(t)\,\overline{\psi\left( \frac{t-b}{a}\right)} \exp\left\{ \frac{i(t^2-b^2)\cot\alpha}{2}\right\} dt\eqno(2.10)$$

\parindent=0mm\vspace{.1in}
where $\psi_{\alpha,a,b}(t)\in L^2(\mathbb R)$ is given by (2.8).

\parindent=0mm\vspace{.2in}
The FrWT (2.10) deals generally with continuous functions, i.e. functions which are defined at all values of the time $t$. However, in many applications, especially in signal processing, data are represented by a finite number of values, so it is important and often useful to consider the discrete version of the continuous FrWT  (2.10). From a mathematical point of view, the continuous parameters $a$ and $b$ in (2.8) can be converted into a discrete one by assuming that $a$ and $b$ take only integral values. For a good discritization of the  wavelets, we choose $a=a_{0}^{-j}$ and $b=kb_{0}a_{0}^{-j}$, where $a_{0}$ and $b_{0}$  are fixed positive constants. Hence, the discritized wavelet family is defined as

$$\psi_{\alpha,j,k}(t)=a_{0}^{j/2}\,\psi\left(a_{0}^{j}t-kb_{0}\right)\exp\left\{-i~\frac{t^2-\big(kb_{0}a_{0}^{-j}\big)^2}{2}\cot\alpha\right\} \eqno(2.11)$$

\parindent=0mm\vspace{.1in}
where the integers $j$ and $k$ are the controlling factors for the dilation and translation, respectively and are contained in a set of integers. For computational efficiency, the discrete wavelet parameters $a_{0}=2$ and $b_{0}=1$ are commonly used so that equation (2.11) becomes
$${\mathscr F}_{\psi}(j,k):=\left\{\psi_{\alpha, j,k}(t)=2^{j/2}\,\psi\left(2^{j}t-k\right)e^{-i~\frac{t^2-(k2^{-j})^2}{2}\cot\alpha}, j,k\in\mathbb Z\right\} . \eqno(2.12)$$

\parindent=0mm\vspace{.1in}
The fractional wavelet system ${\mathscr F}_{\psi}(j,k)$ is called a {\it fractional wavelet frame}, if there exist positive constants $A$ and
$B$ such that

$$A\big\|f \big\|^2_{2} \le \sum_{j\in\mathbb Z}\sum_{k\in \mathbb Z} \left|\big\langle f, \psi_{\alpha, j, k}\big\rangle\right|^2 \le B \big\|f\big\|^2_{2},\eqno(2.13)$$

\parindent=0mm \vspace{.1in}
holds for every $f\in  L^2(\mathbb R)$, and we call the optimal constants $A$ and $B$ the lower frame bound and the upper frame bound, respectively. A {\it tight fractional wavelet frame} refers to the case when $A = B$, and a Parseval  frame refers to the case when $A = B = 1$. On the other hand if only the right hand side of the above double inequality holds, then we say ${\mathscr F}_{\psi}(j,k)$ a {\it Bessel system}.

\parindent=0mm \vspace{.2in}
\textbf{3. Fractional  Multiresolution Analysis on $\mathbb R$.}

\parindent=0mm \vspace{.1in}
We first define a fractional multiresolution analysis on $\mathbb R$ as follows (see \cite{2}):
 
\parindent=0mm \vspace{.1in}
\textbf{Definition 3.1.} A fractional multiresolution analysis is defined as a sequence of closed subspaces $\{V_k^\alpha\} \in L^2(\mathbb R)$ such that

\parindent=0mm \vspace{.1in}
(a) $V_k^\alpha \subseteq V_{k+1}^\alpha,~~k \in \mathbb Z;$

\parindent=0mm \vspace{.1in}
(b) $\bigcup_{k \in \mathbb Z} V_k^\alpha$ is dense in $L^2(\mathbb R)$;

\parindent=0mm \vspace{.1in}
(c) $\bigcap_{k \in \mathbb Z} V_k^\alpha = \{0\}$;

\parindent=0mm \vspace{.1in}
(d) $f(t) \in V_k^\alpha$ if and only if $f(2t)e^{\frac{j}{2}[(2t)^2 - t^2]\cot\alpha} \in V_{k+1}^\alpha,~~k \in \mathbb Z;$

\parindent=0mm \vspace{.1in}
(e) there is a function $ \phi \in V_0^\alpha $, called  {\it fractional scaling function} such that $\{\phi_{\alpha,0,n} = \phi(t-n)e^{-j(tn+n^2)\cot\alpha} : n \in \mathbb Z\}$ is an orthonormal basis of subspace $V_0^\alpha$.

\parindent=8mm \vspace{.1in}
In the above definition, if we assume that the set of functions $\{\phi_{\alpha,0,n} : n\in \mathbb Z\}$ form a Reisz basis of $V_0^\alpha$, then $\phi(t)$ generates a generalized fractional MRA $\{V_m^\alpha\}$ of $L^2(\mathbb R)$, then
$$\phi_{\alpha, m,n}(t) = 2^{\frac{m}{2}}\phi(2^m t-n)e^{\frac{-j}{2}[t^2-(2^{-m}n)^2 -(2^mt-n)^2]\cot\alpha}$$
is the orthonormal basis of $\{V_m^\alpha\}$.

\parindent=8mm \vspace{.1in}
Given a fractional MRA $\{V_j^\alpha : j\in\mathbb Z\}$, we define another sequence $\{W_j^\alpha : j\in\mathbb Z\}$ of closed subspaces of $L^2(\mathbb R)$ by

$$W_j^\alpha = V_{j+1}^\alpha \ominus V_j^\alpha.$$
These subspaces also satisfy 
$$ f\in W_j^\alpha~~\textit{if and only if}~~f(2t)e^{\frac{j}{2}(2t)^2 -t^2]\cot\alpha} \in W_{j+1}^\alpha, ~j\in \mathbb Z.\eqno(3.1)$$
Moreover, they are mutually orthogonal, and we have the following orthogonal decompositions:
\begin{align*}
L^2(\mathbb R) & =\bigoplus_{j\in \mathbb Z}W_j^\alpha\tag{3.2}\\\
&=V_0^\alpha\oplus\left(\bigoplus_{j\ge 0}W_j^\alpha\right)\tag{3.3}
\end{align*}
It should be noted that (3.2) means that the orthonormal basis for $L^2(\mathbb R)$ can be constructed by finding out an orthonormal basis for the subspace $W_j^\alpha$.

\parindent=8mm \vspace{.1in}

The following lemma is very useful in establishing various results and can be found in \cite{dzw}.

\parindent=0mm \vspace{.1in}
\textbf{Lemma 3.2.} The system $\left\{f(t-n)e^{-j(tn+n^2)\cot\alpha} : n\in \mathbb Z\right\}$ of functions is an orthonormal system in $L^2(\mathbb R)$ if and only if
$$\displaystyle\sum_{k\in \mathbb Z}\left|\mathcal{F}_\alpha\{f\}\left(u+2k\pi\sin\alpha\right)\right|^2 = \dfrac{1}{\sin\alpha}.$$

\parindent=8mm \vspace{.1in}
Let $\left\{V_j^\alpha :j \in \mathbb Z\right\}$ be a fractional MRA of $L^2(\mathbb R)$. Since $\phi_{\alpha,0,0}(t) \in V_0^\alpha \subseteq V_1^\alpha,$ and $\{\phi_{\alpha,1,n}(t) : n \in \mathbb Z\}$ is orthonormal basis of $V_1^\alpha$, there must exist coefficient $\{c[k]\}_{k\in \mathbb Z}$ such that
$$\phi_{\alpha,0,0}(t) = \sum_{k \in \mathbb Z} c[k]\phi_{\alpha,1,n}(t)$$
which can be simplified as
$$ \phi(t) = \sum_{n \in \mathbb Z} h[n] \sqrt{2}\phi(2t-n)e^{\frac{-j}{2}[t^2-(\frac{n}{2})^2 - (2t-n)^2]\cot\alpha} \eqno(3.4)$$
and the coefficient can be solved as 
$$ h[n] = \sqrt{2}\int_{-\infty}^{\infty} \phi(t) \phi^{*}(2t-n)e^{\frac{j}{2}[t^2-(\frac{n}{2})^2 - (2t-n)^2]\cot\alpha}dt $$
By taking the FrFT on both sides of Eq. (3.4), we have
\begin{align*}
\Theta_\alpha(u) & =\sum_{n \in \mathbb Z} h[n]\sqrt{2}\mathcal{A}_\alpha\int_{-\infty}^{\infty}\phi(2t-n)e^{\frac{j}{2}[t^2-(\frac{n}{2})^2 - (2t-n)^2]\cot\alpha -jtu\csc\alpha}dt \\\
&=\dfrac{1}{\sqrt{2}}e^{\frac{3ju^2}{8}\cot\alpha}\sum_{n \in \mathbb Z} h[n]e^{\frac{jn^2}{8}\cot\alpha -\frac{jnu}{2}\csc\alpha}\\\
&\qquad\qquad\qquad\times\mathcal{A}_\alpha\int_{-\infty}^{\infty}\phi(2t-n)e^{\frac{j}{2}[t^2-(\frac{n}{2})^2 - (2t-n)^2]\cot\alpha -j(2t-n)u\csc\alpha}d(2t-n) \\\
&=\dfrac{1}{\sqrt{2}}e^{\frac{3ju^2}{8}\cot\alpha}\sum_{n \in \mathbb Z} h[n]e^{\frac{jn^2}{8}\cot\alpha -\frac{jnu}{2}\csc\alpha}\Theta_\alpha\left(\frac{u}{2}\right)\\\
&=\dfrac{1}{\sqrt{2}}e^{\frac{3ju^2}{8}\cot\alpha}C_\alpha\left(\frac{u}{2}\right)\Theta_\alpha\left(\frac{u}{2}\right),\tag{3.5}
\end{align*}
where
$$C_\alpha(u) =\sum_{n \in \mathbb Z} h[n]e^{\frac{jn^2}{8}\cot\alpha -\frac{jnu}{2}\csc\alpha}.$$

By defining
\begin{align*}
\Lambda_\alpha(u) &=\dfrac{1}{\sqrt{2}}e^{\frac{3ju^2}{2}\cot\alpha}C_\alpha(u) \\\\
& =\dfrac{1}{\sqrt{2}}\sum_{n \in \mathbb Z} f[n] \mathcal{A}_\alpha e^{\frac{jn^2}{2}\cot\alpha - jnu\csc\alpha}.
\end{align*}
Eq.(3.5) can be written as 
$$\Theta_\alpha(u) = \Lambda_\alpha\left(\frac{u}{2}\right)\Theta_\alpha\left(\frac{u}{2}\right).\eqno(3.6)$$

It is to be noted that $\Lambda_\alpha(u) $ is  a $2k\pi\sin\alpha$-periodic function since we have
\begin{align*}
\Lambda_\alpha(u+2k\pi\sin\alpha) &= \dfrac{1}{\sqrt{2}}\sum_{n \in \mathbb Z} f[n] \mathcal{A}_\alpha e^{\frac{jn^2}{2}\cot\alpha - jn(u +2k\pi\sin\alpha)\csc\alpha} \\ \\
&=\dfrac{1}{\sqrt{2}}\sum_{n \in \mathbb Z} f[n] \mathcal{A}_\alpha e^{\frac{jn^2}{2}\cot\alpha - jnu\csc\alpha}\\ \\
&= \Lambda_\alpha(u).
\end{align*}
In some of the results in this paper we only need that function of translates is a Reisz basis of $V_0^\alpha,$ which is weaker than being an orthonormal basis. Let $\mathcal{H}$ be a closed subspace of $L^2(\mathbb R)$. A system $\{f_k :k \in \mathbb Z\}$ of functions in $L^2(\mathbb R)$ is said to be a $\textit{Reisz basis}$ of $\mathcal{H}$ if for any $ f \in \mathcal{H}$, there is a sequence $\{a[k] : k \in \mathbb Z\} \in \ell^2(\mathbb Z)$ such that
$$f = \sum_{k \in \mathbb Z}a[k] f_k~~~\textit{with convergence in}~~L^2(\mathbb R),$$ 
and  
$$C_1\sum_{k\in \mathbb Z}\left|a[k]\right|^2 \le \left\| \sum_{k \in \mathbb Z}a[k] f_k\right\|_2^2 \le C_2\sum_{k\in \mathbb Z}\left|a[k]\right|^2,$$  
where the constants $C_1, C_2$ satisfy $ 0 < C_1 \le C_2 <\infty$ and are independent of $f$.

\parindent=0mm \vspace{.1in}
\textbf{Lemma 3.3.} Let $\phi \in L^2(\mathbb R)$ be such that $\{\phi_{\alpha,0,n} = \phi(t-n)e^{-j(tn+n^2)\cot\alpha} : n \in \mathbb Z\}$ forms a Reisz basis of its closed linear span with constants $C_1$ and $C_2$. The for a.e. $u \in [0,2\pi\sin\alpha]$,
$$C_1 \le \sum_{k\in \mathbb Z} \left|\Theta_\alpha\left(u+2k\pi\sin\alpha\right)\right|^2 \le C_2.\eqno(3.7)$$

\textbf{Proof.} By hypothesis, we have 
$$ C_1 \sum_{k \in \mathbb Z} |a[k]|^2 \le \left\| \sum_{k \in \mathbb Z}a[k] \phi(t-n)e^{-j(tn+n^2)\cot\alpha}\right\|_2^2 \le C_2\sum_{k\in \mathbb Z}\left|a[k]\right|^2,\eqno(3.8)$$
where $C_1$ and $C_2$ satisfy $ 0 <C_1 \le C_2 < \infty$ and they are independent of the sequence $\{a[k]\} \in \ell ^2 (\mathbb Z).$

\parindent=0mm \vspace{.1in}
Let us consider the set 
$$\Omega = \left\{u \in \mathbb [0,2\pi\sin\alpha] :\sum_{k\in \mathbb Z} \left|\Theta_\alpha\left(u+2k\pi\sin\alpha\right)\right|^2 > \gamma\right\}$$
Here we assume that $\Omega$ has positive measure. We will show that $\gamma \le C_2.$ Consider a sequence $\left\{c[k]\right\} \in \ell^2(\mathbb Z)$ such that 
$$\chi_{\Omega}(u)=\sum_{k \in \mathbb Z} a[k]e^{-2j\pi u}~~\textit{for a.e.}~ u\in [0,2\pi\sin\alpha].$$
Then it implies that
\begin{align*}
\left\| \sum_{k \in \mathbb Z}a[k] \phi(t-n)e^{-j(tn+n^2)\cot\alpha}\right\|_2^2& = \int_{-\infty}^{\infty} \left|\sum_{k \in \mathbb Z} a[k]e^{-2j\pi u}\right|^2\,|\Theta_\alpha(u)|^2\,du\\\\
&= \int_{0}^{2\pi\sin\alpha} \left|\sum_{k \in \mathbb Z} a[k]e^{-2j\pi u}\right|^2\,\sum_{\ell \in \mathbb Z}|\Theta_\alpha(u+2\pi\ell\sin\alpha)|^2\,du\\\\
&=\int_{\chi_{\Omega}}\sum_{\ell \in \mathbb Z}|\Theta_\alpha(u+2\pi\ell\sin\alpha)|^2\,du\\\
& \ge \int_{\chi_{\Omega}} \gamma du\\
& =\gamma |\Omega|
\end{align*}
By Parseval's identity, we have $\sum_{k \in \mathbb Z} |a[k]|^2 = |\Omega|.$ Hence, 
$$\left\| \sum_{k \in \mathbb Z}a[k] \phi(t-n)e^{-j(tn+n^2)\cot\alpha}\right\|_2^2 \ge \gamma|\Omega| = \gamma\sum_{k \in \mathbb Z} |a[k]|^2.$$

On comparing with (3.8), it is clear that $\gamma \le C_2,$ as required. Hence the set 
$$\left\{u \in \mathbb [0,2\pi\sin\alpha] :\sum_{k\in \mathbb Z} \left|\Theta_\alpha\left(u+2k\pi\sin\alpha\right)\right|^2 > C_2\right\}$$ ha s measure zero. Therefore,
$$\sum_{k\in \mathbb Z} \left|\Theta_\alpha\left(u+2k\pi\sin\alpha\right)\right|^2 \le  C_2~~~\textit{for a.e.}~~ u \in [0,2\pi\sin\alpha].$$
Similarly, considering the set
$$\Xi = \left\{u \in \mathbb [0,2\pi\sin\alpha] :\sum_{k\in \mathbb Z} \left|\Theta_\alpha\left(u+2k\pi\sin\alpha\right)\right|^2 < \gamma\right\}$$
we get the left hand inequality of (3.7).  $\square$

\parindent=0mm \vspace{.1in}
\textbf{Lemma 3.4.} Let $\phi \in L^2(\mathbb R)$ such that the collection  $\{\phi_{\alpha,0,n} (t) : n \in \mathbb Z\}$ is a Reisz basis of the space
$$V_0^\alpha = \left\{\sum_{n \in\mathbb Z} c[n] \phi_{\alpha,0,n}(t) : c[n] \in \ell^2(\mathbb Z)\right\}$$
of $L^2(\mathbb R)$ if and only if there exists positive constants $ A, B$ such that for all $ u\in I = [0, 2\pi\sin\alpha]$, we have

$$ A \le \mathcal{G}^2_{\alpha, \phi}(u) \le B \eqno(3.9)$$

where

$$\mathcal{G}_{\alpha, \phi }(u) = \sqrt{2\pi\sin\alpha\sum_{k \in \mathbb Z}\left|\Theta_\alpha(u+2k\pi\sin\alpha)\right|^2}.\eqno(3.10)$$

\parindent=0mm \vspace{.1in}
\textbf{Proof.} For any $f(t) \in V_0^\alpha$, we have
$$ f(t) = \sum_{n\in \mathbb Z} c[n] \phi_{\alpha,0,n}(t)\eqno(3.11)$$
where $c[n]\in\ell^2(\mathbb Z)$.

\parindent=0mm \vspace{.1in}
On taking FrFT on both sides of $(3.11)$, we obtain
$$ \mathcal{F}_\alpha\{f(t)\}(u) = \sqrt{2\pi}\, \widetilde{c}_\alpha(u)\Theta_\alpha(u)\eqno(3.12)$$
where $\widetilde{c}_\alpha(u)$ denotes the discrete FrFT of $c[n]$. By using Parseval formula of the FrFT, we have
\begin{align*}
\|f(t)\|_{L^2(\mathbb R)}^2 & = \|\mathcal{F}_\alpha\{f(t)\}(u)\|_{L^2(\mathbb R)}^2\\\\
&= \int_{-\infty}^{\infty} 2\pi |\widetilde{c}_\alpha(u)|^2\,|\Theta_\alpha(u)|^2\,du\\\\
&= \sum_{k \in \mathbb Z}\int_{0}^{2\pi\sin\alpha}2\pi |\widetilde{c}_\alpha(u+2k\pi\sin\alpha)|^2\,|\Theta_\alpha(u+2k\pi\sin\alpha)|^2\,du\\\\
& =\int_{0}^{2\pi\sin\alpha}|\widetilde{c}_\alpha(u)|^2\mathcal{G}^2(\alpha, \phi,u)\,du.\tag{3.13}
\end{align*}
Further, Parsevals formula for discrete FrFT yields
$$\|c[n]\|_{\ell^2(\mathbb Z)}^2 = \sum_{n \in \mathbb Z}|c[n]|^2 = \int_{0}^{2\pi\sin\alpha}|\widetilde{c}_\alpha(u)|^2\,du\eqno(3.14)$$
Now, Eqns. $(3.9), (3.13)$  and $ (3.14)$ yields
$$ A \|c[n]\|_{\ell^2(\mathbb Z)}^2 \le \left\|\sum_{n \in \mathbb Z}c[n]\phi_{\alpha,0,n}(t)\right\|^2 \le B \|c[n]\|_{\ell^2(\mathbb Z)}^2.\eqno(3.15)$$
It follows from $(3.15)$  that $\{\phi_{\alpha,0,n}(t) : n \in \mathbb Z\}$ is a Reisz basis for $V_0^\alpha$. In particular $\{\phi_{\alpha,0,n}(t) : n \in \mathbb Z\}$ is an orthonormal basis for $V_0^\alpha$ if and only if $ A = B = 1.\square$

\parindent=0mm \vspace{.1in}
Now the following results shows that if the  fractional translates of a function  form a REisz basis for the spanned subspace, then there exists another function whose fractional translates form an orthonormal basis for the same subspace.

\parindent=0mm \vspace{.1in}
{\bf{Theorem 3.5.}} Suppose  that $\{ \phi(t-n)e^{-j(tn+n^2)\cot\alpha} : n \in \mathbb Z\}$ forms a Reisz basis of its closed linear span $V_0^\alpha$. Then there is a function $\varphi$ such that $\{\varphi_{\alpha,0,n} = \varphi(t-n)e^{-j(tn+n^2)\cot\alpha} : n \in \mathbb Z\}$ forms an orthonormal basis for $V_0^\alpha.$ 

\parindent=0mm \vspace{.1in}
\textbf{Proof.} By Lemma 3.4., if $A = B = 1,~\{\phi_{\alpha, 0, n}(t) : n\in \mathbb Z \}$ is an orthonormal basis for $V_0^\alpha$ so that $\mathcal{G}_{\alpha,\phi} (u) = 1.$ Define $\varphi$ so that
$$\mathcal{F}_\alpha\{\phi(t)\}(u) = \dfrac{\Theta_\alpha(u)}{\sqrt{2\pi\sin\alpha\sum_{k \in \mathbb Z}\left|\Theta_\alpha(u+2k\pi\sin\alpha)\right|^2}}\eqno(3.16).$$

\parindent=0mm \vspace{.1in}
It follows from (3.7) that $ \mathcal{F}_\alpha\{\phi\} \in L^2(\mathbb R)$. Hence $\varphi$ also belongs to $L^2(\mathbb R)$.It is clear from (3.16) that $\varphi(t) = \phi(t)$. Hence $\{\varphi_{\alpha,0,n} : n \in \mathbb Z\}$ is an orthonormal basis for $V_0^\alpha$ if and only if $A = B = 1.$ Otherwise,  let $\varphi \in L^2(\mathbb R)$, then there exists a sequence $\{c[n]\}_{n \in \mathbb Z} \in \ell^2(\mathbb Z)$ satisfies
\begin{align*}
\varphi(t)&= \sum_{n \in \mathbb Z} c[n] \phi_{\alpha,0,n}(t)\\\
&=\sum_{n \in \mathbb Z} c[n]\phi(t-n)e^{-j(tn+n^2)\cot\alpha}
\end{align*}
Taking FrFT on both sides, we obtain
$$\mathcal{F}_\alpha\{\varphi(t)\}(u) = \sqrt{2\pi}\sin\alpha \Gamma_\alpha(u) \Theta_\alpha(u). \eqno(3.17)$$

where $\Gamma_\alpha(u)$ denotes the DTFT of $c[n]e^{jn^2\cot\alpha}$ and is $2\pi\sin\alpha$ periodic. Then by utilizing $(3.16)$ and $(3.17)$, we obtain
$$\Gamma_\alpha(u) = \dfrac{1}{\sqrt{2\pi}\sin\alpha\mathcal{G}_{\alpha,\phi}(u)}\eqno(3.18)$$
Meanwhile, since $\Gamma_\alpha(u)$ is $2\pi\sin\alpha$ periodic, applying (3.17) yields
$$\mathcal{F}_\alpha\{\varphi\}(u+2k\pi\sin\alpha) = \sqrt{2\pi}\sin\alpha\Gamma_\alpha(u)\Theta_\alpha(u+2k\pi\sin\alpha)$$
so that 
$$\left|\mathcal{F}_\alpha\{\varphi\}(u+2k\pi\sin\alpha)\right|^2 = 2\pi\sin\alpha\left|\Gamma_\alpha(u)\right|^2\,\left|\mathcal{F}_\alpha\{\varphi\}(u+2k\pi\sin\alpha)\right|^2.\eqno(3.19)$$
Taking a sum for all $k$ on both sides of (3.19) yields
$$\mathcal{G}_{\alpha,\varphi}^2(u)= 2\pi\sin\alpha\left|\Gamma_\alpha(u)\right|^2\,\mathcal{G}_{\alpha,\phi}^2(u)\eqno(3.20)$$
where $\mathcal{G}_{\alpha,\varphi}^2(u) = \sqrt{2\pi\sin\alpha\sum_{k \in \mathbb Z}\left|\mathcal{F}_\alpha\{\varphi\}(u+2k\pi\sin\alpha)\right|^2}.$ Inserting (3.18) into (3.20) yields $\mathcal{G}_{\alpha,\varphi}(u) = 1.$ Then it follows from Lemma 3.4 that $\{\varphi_{\alpha,0,n}(t)\}_{n \in \mathbb Z}$ forms an orthonormal basis for $V_0^\alpha$. This completes the theorem.$\square$

\parindent=0mm \vspace{.1in}
\textbf{4. Union density and Intersection triviality conditions}

\parindent=0mm \vspace{.1in}
\textbf{Theorem 4.1.} Let $\{V_j^\alpha : j \in \mathbb Z\}$ be a sequence of closed subspaces of $L^2(\mathbb R)$ satisfying conditions (a), (d) and (e) of Definition 3.1.Then,

\parindent=0mm \vspace{.1in}
$$\bigcap_{j \in \mathbb Z} V_j^\alpha = \{0\}.$$

\parindent=0mm \vspace{.1in}
This is the case even if, in (e), we only assume that $\{\phi(t-n)e^{-j(tn+n^2)\cot\alpha} : n \in \mathbb Z\}$ is a Reisz basis.

\parindent=0mm \vspace{.1in}
\textbf{Proof.} Suppose that there exists a non- zero $f \in \bigcap_{j \in \mathbb Z}V_j^\alpha,$ we can assume that $\|f\|_2 = 1$. In particular, $f\in V_{-j}^\alpha$ for each $j \in \mathbb Z$, hence if we let $f_j(t) =2^{\frac{j}{2}}f(2^jt)$ we must have $f_j \in V_0^\alpha$. Further a simple change of variables shows $\|f_j\| = \|f\| =1$. Since we are assuming  that $\{\phi(t-n)e^{-j(tn+n^2)\cot\alpha} : n \in \mathbb Z\}$ is a Reisz basis, we can write 

\parindent=0mm \vspace{.1in}
$$f_j(t) = \sum_{k \in \mathbb Z}a^j[k]\phi_{\alpha,0,k},$$

with the convergence in $L^2(\mathbb R)$, in such a way that
$$A\sum_{k \in \mathbb Z}\left|a^j[k]\right|^2 \le \|f_j\|_2^2 = 1.$$

Taking fractional Fourier transform (FrFT), we can obtain

$$\mathcal{F}_\alpha\{f_j\}(u) =\Lambda_\alpha^j(u)\Theta_\alpha(u),$$
where,

$$\Lambda_\alpha^j(u) = \dfrac{1}{\sqrt{2}}\sum_{k \in \mathbb Z}a^j[k]e^{\frac{jn^2\cot\alpha}{2}-jnu\csc\alpha},$$

is a $2\pi\sin\alpha$-periodic function which belong to $L^2[0,2\pi\sin\alpha]$ with norm $\le \sqrt{\frac{2\pi\sin\alpha}{A}}$. Thus we have,

$$\mathcal{F}_\alpha\{f\}(u) = 2^{\frac{j}{2}} \Lambda_\alpha^j\left(2^ju\right) \Theta_\alpha\left(2^ju\right).$$
Also for $j \ge 1,$ we have
\begin{align*}
\int_{2\pi\sin\alpha}^{4\pi\sin\alpha} \left|\mathcal{F}_\alpha\{f\}(u)\right|\,du &\le 2^{\frac{j}{2}}\left\{\int_{2\pi\sin\alpha}^{4\pi\sin\alpha} \left|\Theta_\alpha(2^j u)\right|^2\,du\right\}^{1/2}\left\{\int_{2\pi\sin\alpha}^{4\pi\sin\alpha} \left|\Lambda_\alpha^j(2^j u)\right|^2\,du\right\}^{1/2}\\\\
& = 2^{\frac{-j}{2}}\left\{\int_{2^{j+1}\pi\sin\alpha}^{2^{j+2}\pi\sin\alpha} \left|\Theta_\alpha(2^j u)\right|^2\,du\right\}^{1/2}\left\{\int_{2^{j+1}\pi\sin\alpha}^{2^{j+2}\pi\sin\alpha} \left|\Lambda_\alpha^j(2^j u)\right|^2\,du\right\}^{1/2}\\\\
& \le \left\{\int_{2^{j+1}\pi\sin\alpha}^{\infty} \left|\Theta_\alpha(2^j u)\right|^2\,du\right\}^{1/2}\left\{\dfrac{1}{2^j}\int_{2^{j+1}\pi\sin\alpha}^{2^{j+2}\pi\sin\alpha} \left|\Lambda_\alpha^j(2^j u)\right|^2\,du\right\}^{1/2}\\\\
& = \left\{\int_{2^{j+1}\pi\sin\alpha}^{\infty} \left|\Theta_\alpha(2^j u)\right|^2\,du\right\}^{1/2}\\\
&\qquad\qquad\qquad\qquad\times\left\{\dfrac{1}{2^j}\sum_{\ell =0}^{2^j -1}\int_{2^{j+1}\pi\sin\alpha+2\ell\pi\sin\alpha}^{2^{j+2}\pi\sin\alpha+2(\ell+1)\pi\sin\alpha} \left|\Lambda_\alpha^j(2^j u)\right|^2\,du\right\}^{1/2}\\\\
& \le \left\{\int_{2^{j+1}\pi\sin\alpha}^{\infty} \left|\Theta_\alpha(2^j u)\right|^2\,du\right\}^{1/2}\sqrt{\dfrac{2\pi\sin\alpha}{A}}.
\end{align*}

\parindent=0mm \vspace{.1in}
By letting $j \rightarrow \infty,$ we obtain $\displaystyle\int_{2\pi\sin\alpha}^{4\pi\sin\alpha} \left|\mathcal{F}_\alpha\{f\}(u)\right|\,du=0$ and we can deduce that 
$$\mathcal{F}_\alpha\{f\}(u) =0 ~\textit{ a.e, on }~[2\pi\sin\alpha, 4\pi\sin\alpha].$$
On applying the same argument to $2^{\frac{\ell}{2}}\mathcal{F}_\alpha\{f\}(2^\ell u),\ell \in \mathbb Z,$ to obtain $\mathcal{F}_\alpha\{f\}(u)=0$ a.e on $2^\ell[2\pi\sin\alpha, 4\pi\sin\alpha], \ell \in \mathbb Z$. Hence $\mathcal{F}_\alpha\{f\}(u) =0$ a.e on $(0,\infty).$ If we apply this argument, with the interval $[-4\pi\sin\alpha,-2\pi\sin\alpha]$ playing the role of $[2\pi\sin\alpha, 4\pi\sin\alpha],$ we obtain $\mathcal{F}_\alpha\{f\}(u) =0$ a.e, on $(-\infty,0)$.~~~~$\square$ 

\parindent=0mm \vspace{.1in}
\textbf{Theorem 4.2.} Let $\{V_j^\alpha : j \in \mathbb Z\}$ be a sequence of a closed subspaces of $L^2(\mathbb R)$ satisfying conditions (a), (d) and (e) of Def. 3.1. Further, assume that the function $\phi$ of condition (e) is such that $\Theta_\alpha$ is continuous at $u=0$. Then the following conditions are equivalent:

\parindent=0mm \vspace{.1in}
(i)~~$\Theta_\alpha(0) \ne 0.$

\parindent=0mm \vspace{.1in}
(ii)~~$ \overline{\bigcup_{j \in \mathbb Z}V_j^\alpha} = L^2(\mathbb R).$

\parindent=0mm \vspace{.1in}
Moreover, when either is the case, we have $\left|\Theta_\alpha(0)\right| =1.$

\parindent=0mm \vspace{.1in}
\textbf{Proof.} Assume that $\Theta_\alpha(0) \ne 0$. we first claim that  $W^\alpha = \overline{\bigcup_{j \in \mathbb Z} V_j^\alpha}$ is invariant under translations. To prove this we first show that $W^\alpha$ is invariant under the dyadic translations $\mathcal{T}_{2^{-\ell}d},~\ell, d \in \mathbb Z.$ Let $f \in W^\alpha$, therefore for given $\epsilon >0,$ there exists $j_0 \in \mathbb Z$ and $h \in V_{j_0}^\alpha$ such that $\|f-h\|_2 < \epsilon.$ From (a) we can deduce that $h \in V_j^\alpha$ for all $j \ge j_0$ and using (d) and (e), we can write

\parindent=0mm \vspace{.1in}
$$h(t) =\sum_{k \in \mathbb Z}c^j[k]\phi(2^j t-k)e^{\frac{-j}{2}[t^2-(2^{-j}k)^2 -(2^jt-k)^2]\cot\alpha}$$

with convergence in $L^2(\mathbb R).$ Hence ,
\begin{align*}
(\mathcal{T}_{2^{-\ell}d} h)(t) &= h(t-2^{-\ell}d)\\\\
&=\sum_{k \in \mathbb Z}c^j[k]\phi\left(2^j (t-2^{-\ell}d)-k\right)e^{\frac{-j}{2}[(t-2^{-\ell}d)^2-(2^{-j}k)^2 -(2^j(t-2^{-\ell}d)-k)^2]\cot\alpha}.
\end{align*}

If $j \ge \ell$ then, $\phi\left(2^j (t-2^{-\ell}d)-k\right)e^{\frac{-j}{2}[(t-2^{-\ell}d)^2-(2^{-j}k)^2 -(2^j(t-2^{-\ell}d)-k)^2]\cot\alpha}$ is an element of $V_j^\alpha$, since $2^{j-\ell}d \in \mathbb Z.$ Since  

$$\left\|\mathcal{T}_{2^{-\ell}d} f - \mathcal{T}_{2^{-\ell}d} h\right\|_2 = \|f-h\|_2 < \epsilon$$

and $\epsilon$ is arbitrarily small, we can conclude that $W^\alpha$ is invariant under dyadic translations. Now, for a general $\tau \in \mathbb R$, we can find integers $d$ and $\ell$ such that $2^{-\ell}d$ is arbitrarily close to $\tau$ , hence we can write

$$\left\|\mathcal{T}_{2^{-\ell}d} f - \mathcal{T}_{\tau} f\right\|_2  < \epsilon$$

\parindent=0mm \vspace{.1in}
and it follows that $W^\alpha$ is invariant under all translations $\mathcal{T}_\tau$.

\parindent=8mm \vspace{.1in}
Further, as $\Theta_\alpha(0) \ne 0$ and $|\Theta|$ is continuous at $0,~~\Theta_\alpha(u) \ne 0 $ on $(-\sigma,\sigma)$ for some $\sigma >0.$ Suppose that there exists $g \in (W^\alpha)^{\perp},$ then g is orthogonal to all $f \in W^\alpha,$ and since $W^\alpha$ is translation invariant, we have

\parindent=0mm \vspace{.1in}
$$\int_{-\infty}^{\infty} f(t+\tau) \overline{g(t)} \,dt = 0$$

\parindent=0mm \vspace{.1in}
for all $\tau \in \mathbb R$ and all $f \in W^\alpha$. This equality and the Plancherel formula implies

\parindent=0mm \vspace{.1in}
$$\int_{-\infty}^{\infty} \mathcal{F}_\alpha\{f\}(u) \overline{\mathcal{F}_\alpha\{g\}(u)}\,du =0.$$

\parindent=0mm \vspace{.1in}
Since $\mathcal{F}_\alpha\{f\}\overline{\mathcal{F}_\alpha\{g\}} \in L^1(\mathbb R)$ this shows that $\mathcal{F}_\alpha\{f\}(u) \overline{\mathcal{F}_\alpha\{g\}(u)} =0$ for a.e. $ u \in \mathbb R$. In particular, letting $f(t) = 2^j\phi(2^j t)e^{\frac{-j}{2}[t^2 -(2^jt)^2]\cot\alpha}$ so that $f \in V_j^\alpha \subset W^\alpha$ and $\mathcal{F}_\alpha\{f\}(u) = \Theta_\alpha\left(2^{-j}u\right).$ Hence $\Theta_\alpha\left(2^{-j}u\right)\overline{\mathcal{F}_\alpha\{g\}(u)} =0 $ for a.e. $u \in \mathbb R$. Since $\Theta_\alpha\left(2^{-j}u\right)\ne 0$ if $ u \in (-2^j\sigma,2^j\sigma),$ we can conclude that $\mathcal{F}_\alpha\{g\}(u)=0$ for a.e $|u| < 2^j\sigma.$ Letting $j \rightarrow \infty,$ we see that $\mathcal{F}_\alpha\{g\}=0$ a.e. and therefore $g =0$. This shows that $\overline{\bigcup_{j \in\mathbb Z} V_j^\alpha } = L^2(\mathbb R).$

\parindent=8mm \vspace{.1in}
Now we assume that $W^\alpha = \overline{\bigcup_{j \in\mathbb Z} V_j^\alpha } = L^2(\mathbb R)$ and let $f$ be such that $\mathcal{F}_\alpha\{f\} = \chi_{[-1,1]},$ then $\|f\|_2^2 =\dfrac{1}{2\pi\sin\alpha}\left\|\mathcal{F}_\alpha\{f\}\right\|_2^2 = \dfrac{1}{\pi\sin\alpha}.$ If  $\mathcal{P}_j$ denotes the orthogonal projection onto $V_j^\alpha$ then we have $\|f -\mathcal{P}_jf\|_2 \rightarrow 0$ as $ j\rightarrow \infty$, due to (a) and our assumption. Thus $\|\mathcal{P}_j\|_2 \rightarrow\|f\|_2$ as $j \rightarrow \infty$. Therefore, if $ \phi_{\alpha,j,k} = 2^{\frac{j}{2}}\phi(2^j t-k)e^{\frac{-j}{2}[t^2-(2^{-j}k)^2 -(2^jt-k)^2]\cot\alpha},$ we have

\parindent=0mm \vspace{.1in}
$$\left\|\mathcal{P}_jf\right\|_2^2 = \left\|\sum_{k \in \mathbb Z}\left\langle f, \phi_{\alpha,j,k}\right\rangle \phi_{\alpha,j,k}\right\|_2^2 \rightarrow\dfrac{1}{\pi\sin\alpha}$$

\parindent=0mm \vspace{.1in}
as $ j \rightarrow \infty,$ since $\{\phi_{\alpha,j,k} : k\in \mathbb Z\}$ is an orthonormal basis of $V_j^\alpha$. By virtue of Plancherel theorem and the fact that $\mathcal{F}_\alpha\{f\} = \chi_{[-1,1]},$ we have
\begin{align*}
\dfrac{1}{(2\pi\sin\alpha)^2}&\sum_{k\in \mathbb Z}\left|\int_{-\infty}^{\infty} \mathcal{F}_\alpha\{f\}(u)\overline{\mathcal{F}_\alpha\{\phi_{\alpha,j,-k}}(u)\}\,du\right|^2\\
&\qquad =\dfrac{1}{(2\pi\sin\alpha)^2}\sum_{k\in \mathbb Z}\left|\int_{-\infty}^{\infty} \mathcal{F}_\alpha\{f\}(u)2^{\frac{-j}{2}} e^{-j2^{-j}ku\csc\alpha}\overline{\Theta_\alpha\left(2^{-j}u\right)}\,du\right|^2\\ \\
&\qquad=2^j \sum_{k\in \mathbb Z}\left|\dfrac{1}{2\pi\sin\alpha}\int_{-2^{-j}\pi\sin\alpha}^{2^{-j}\pi\sin\alpha}\overline{\Theta_\alpha\left(\xi\right)}e^{-jk\xi\csc\alpha}\,d\xi\right|^2.
\end{align*}

\parindent=0mm \vspace{.1in}
For large enough $j, [-2^{-j}\pi\sin\alpha,2^{-j}\pi\sin\alpha] \subset [-\pi\sin\alpha,\pi\sin\alpha]$ and the last expression is $2^j$ times the sum of squares of the absolute values of the Fourier coefficients of the function $\chi_{[-2^{-j},2^{-j}]}\overline{\Theta_\alpha},$  thus by invoking Plancherel formula for fourier series,

\parindent=0mm \vspace{.1in}
$$ \dfrac{2^j}{2\pi\sin\alpha}\int_{-2^{-j}\pi\sin\alpha}^{2^{-j}\pi\sin\alpha}\left|\Theta_\alpha(\xi)\right|^2\,d\xi \rightarrow \dfrac{1}{\pi\sin\alpha}$$

\parindent=0mm \vspace{.1in}
as $ j\rightarrow \infty.$ But, by the continuity of $|\Theta_\alpha|$ at $0$, the last integral expression tends to $\dfrac{1}{\pi\sin\alpha}\left|\Theta_\alpha(0)\right|^2.$ Therefore, $\left|\Theta_\alpha(0)\right| = 1 \ne 0.\square$

\parindent=0mm \vspace{.2in}

{\bf{5. Characterization of Fractional Scaling Function}}

\parindent=0mm \vspace{.1in}
In this section we will characterize those functions that are scaling functions for fractional MRA of $L^2(\mathbb R)$. For that we should clarify what we mean by scaling function for  fractional MRA.

\vspace{.1in}\parindent=0mm
For a given function $\phi \in L^2(\mathbb R)$, we define the closed subspaces $\{V_j : j \in \mathbb Z\}$ of $L^2(\mathbb R)$ as follows

\parindent=0mm \vspace{.1in}
$$V_0^\alpha = \overline{span}\{\phi(t-n)e^{-j(tn+n^2)\cot\alpha} : k \in \mathbb Z\}$$
and
$$V_j^\alpha = \left\{ f : f(2^{-j}t)e^{\frac{-j}{2}[(2^{-j}t)^2-t^2]\cot\alpha}\in V_0^\alpha\right\}~~~\textit{for}~~j \in \mathbb Z\backslash\{0\}.$$
we say that $\phi \in L^2(\mathbb R)$ is a scaling function for a fractional multiresolution analysis of $L^2(\mathbb R)$ if the sequence of fractional closed subspaces $\{V_j^\alpha : j \in \mathbb Z\}$ as defined above constitutes a fractional multiresolution analysis of $L^2(\mathbb R)$.

\parindent=0mm \vspace{.1in}
\textbf{Theorem 5.1.} A function $\phi \in L^2(\mathbb R)$ is a fractional scaling function for a fractional MRA if and only if 

$$\sum_{k\in \mathbb Z} \left|\Theta_\alpha(u+2k\pi\sin\alpha)\right|^2 = \dfrac{1}{\sin\alpha}~~~\textit{for a.e.}~~u \in [0,2\pi\sin\alpha];\eqno(5.1)$$

$$\lim_{j\rightarrow \infty}\left|\Theta_\alpha(2^{-j}u)\right| = 1~~\textit{for a.e.}~~ u \in \mathbb R; \eqno(5.2)$$

and there exists a $2\pi\sin\alpha$-periodic function $\Lambda_\alpha$ such that

$$\Theta_\alpha(2u) = \Lambda_\alpha(u) \Theta_\alpha(u)~~\textit{for a.e.}~~u \in \mathbb R\eqno(5.3)$$

\parindent=0mm \vspace{.1in}
\textbf{Proof.} Suppose that $\phi$ is a scaling function for fractional MRA. Then $\{\phi_{\alpha,0,n} = \phi(t-n)e^{-j(tn+n^2)\cot\alpha} : n \in \mathbb Z\}$ is an orthonormal system in $L^2(\mathbb R)$ which is equivalent to $(5.1)$ by Theorem  3.2. Also $(5.3)$ follows from $(3.6)$ . Now we proceed to prove $(5.2)$ as follows. Since $\{V_j^\alpha :j \in \mathbb Z\}$ is an MRA for $L^2(\mathbb R)$, we have $\overline{\bigcup_{j\in\mathbb Z}V_j^\alpha} = L^2(\mathbb R).$ Following the second part of the proof of the Theorem (), we have

\parindent=0mm \vspace{.0in}
$$\lim_{j\rightarrow \infty}\dfrac{1}{2\pi\sin\alpha}\int_{0}^{2\pi\sin\alpha}\left|\Theta_\alpha\left(2^{-j}u\right)\right|\,du = 1.$$

\parindent=0mm \vspace{.1in}
Using $(5.3)$ in (5.1) , we get
\begin{align*}
\dfrac{1}{\sin\alpha} & =\sum_{k \in \mathbb Z}\left|\Theta_\alpha\left(2u +2k\pi\sin\alpha\right)\right|^2\\\\
&=\sum_{k \in \mathbb Z}\left|\Lambda_\alpha\left(u +k\pi\sin\alpha\right)\right|^2\left|\Theta_\alpha\left(u +k\pi\sin\alpha\right)\right|^2\\\\
&=\left|\Lambda_\alpha(u)\right|^2\sum_{\ell \in \mathbb Z}\left|\Theta_\alpha\left(u +2\ell\pi\sin\alpha\right)\right|^2+\left|\Lambda_\alpha\left(u +k\pi\sin\alpha\right)\right|^2\sum_{\ell \in \mathbb Z}\left|\Theta_\alpha\left(u +(2\ell+1)\pi\sin\alpha\right)\right|^2.
\end{align*}

Thus for all $u \in \mathbb R$, we obtain
$$\left|\Lambda_\alpha(u)\right|^2 +\left|\Lambda_\alpha\left(u+\pi\sin\alpha\right)\right|^2=1.$$

In particular, this shows that
 $$\left|\Lambda_\alpha(u)\right|\le 1~~~~\textit{ for a.e.,}~~~ u \in \mathbb R$$.

This inequality and (3.4) shows that $\left|\Theta_\alpha(2^{-j}u)\right|$ is non decreasing for almost every $u \in \mathbb R$ as $j\rightarrow \infty$.

\parindent=0mm \vspace{.1in}
Letting
$$g(u) = \lim_{j\rightarrow \infty}\left|\Theta_\alpha(2^{-j}u)\right|,$$
Since by Theorem 3.2, we have $\left|\Theta_\alpha(u)\right| \le 1$ a.e., by Lebesgue dominated convergence theorem it now follows that

$$\dfrac{1}{2\pi\sin\alpha}\int_{0}^{2\pi\sin\alpha}g(u)\,du = 1,$$

and (5.2) then follows as for all $u \in \mathbb R,~~ 0 \le g(u) \le 1.$

\parindent=8mm \vspace{.1in}
We now proceed to prove the converse. Assume that (5.1),(5.2) and (5.3) are satisfied. The orthonormality of $\{\phi(t-n)e^{-j(tn+n^2)\cot\alpha} : k \in \mathbb Z\}$ is equivalent to (5.1), as observed earlier. This fact along with the definition of $V_0^\alpha$ gives us (e) of the definition of fractional MRA.

\parindent=8mm \vspace{.1in}
The definition of the subspaces $V_j^\alpha$ also show that $f \in V_j^\alpha$ holds if and only if $f(2t)e^{\frac{j}{2}[(2t)^2 - t^2]\cot\alpha} \in V_{j+1}^\alpha,$ which is (d) of the definition of the fractional MRA. Now, for each $j \in \mathbb Z,$ we claim

$$ V_j^\alpha = \left\{ f :\hat{f}_\alpha(2^j u) = \mu_j(u)\Theta_\alpha(u)~~\textit{for some}~~2\pi\sin\alpha-periodic ~function ~\mu_j\in [0,2\pi\sin\alpha]\right\}.\eqno(5.4)$$
This claim is established by expressing $f(2^{-j}t)e^{\frac{j}{2}(2^{-j}t)^2-t^2]\cot\alpha} \in V_0^\alpha$ as a linear combination of $\{\phi(t-n)e^{-j(tn+n^2)\cot\alpha} : k \in \mathbb Z\}$ and then taking Fourier transforms.

\parindent=8mm \vspace{.1in}
To prove the inclusion $V_j^\alpha \subset V_{j+1}^\alpha,$ it is enough to show that $V_0^\alpha \subset V_0^\alpha$. By (5.4), given $ f\in V_0^\alpha,$ there is an $2\pi\sin\alpha-$ periodic function $\mu_0 \in L^2[0,2\pi\sin\alpha]$ such that 

$$\hat{f}_\alpha(2u)  = \mu_0(2u)\Theta_\alpha(2u).$$
Thus, using (5.3), we get

$$\hat{f}_\alpha(2u) =\mu_0(2u)\Lambda_\alpha(u)\Theta_\alpha(u).$$
It is clear that the function $\mu_0(2u)\Lambda_\alpha(u)$ is $2\pi\sin\alpha$- periodic function. Now

$$\int_{0}^{2\pi\sin\alpha}\left|\mu_0(2u)\right|^2\left|\Lambda_\alpha(u)\right|^2\,du \le \int_{0}^{2\pi\sin\alpha}\left|\mu_0(2u)\right|^2 < \infty,$$

\parindent=0mm \vspace{.1in}
as $\left|\Lambda_\alpha(u)\right| \le 1 $ for a.e. $ u \in [0,2\pi\sin\alpha]$. Hence, the function $\mu_0(2u)\Lambda_\alpha(u)$ belongs to $L^2[0,2\pi\sin\alpha]$. Again by (5.4), $f \in V_1^\alpha.$

\parindent=8mm \vspace{.1in}
We have already seen in Theorem 4.1 that property (c) in the definition of fractional MRA follows from (a), (d) and (e). Now it remains to prove only one property, i.e., we have to show that

\parindent=0mm \vspace{.1in}
$$\overline{\bigcup _{j \in \mathbb Z}V_j^\alpha} = L^2(\mathbb R).$$

\parindent=0mm \vspace{.1in}
 Let $\mathcal{P}_j$ be the projection on $V_j^\alpha$. It suffices to show that

\parindent=0mm \vspace{.1in}
$$\left\|\mathcal{P}_jf-f\right\|_2^2 = \left\|f\right\|_2^2 -\left\|\mathcal{P}_jf\right\|_2^2 \rightarrow 0~~~\textit{as}~~j \rightarrow \infty.$$

\parindent=0mm \vspace{.1in}
Furthermore, we assume that our $f\in L^2(\mathbb R)$ is such that $\hat{f}_\alpha$ has compact support. Since
 $\left\{2^{\frac{j}{2}}\phi(2^j t-k)e^{\frac{-j}{2}[t^2-(2^{-j}k)^2 -(2^jt-k)^2]\cot\alpha} : k \in \mathbb Z\right\}$ is an orthonormal basis of $V_j^\alpha$ and for large positive $ j, \hat{f}_\alpha(2^ju)$ has support in $[-\pi\sin\alpha.\pi\sin\alpha],$

\parindent=0mm \vspace{.1in}
\begin{align*}
\left\|\mathcal{P}_jf\right\|_2^2 &= 2^{-j}\sum_{k \in \mathbb Z}\left|\int_{-\infty}^{\infty} f(t)\overline{\phi(2^j t-k)e^{\frac{-j}{2}[t^2-(2^{-j}k)^2 -(2^jt-k)^2]\cot\alpha}}2^j\,dt\right|^2\\\\
&=2^{-j}\sum_{k \in \mathbb Z}\left|\int_{-\infty}^{\infty} f(2^{-j}t)\overline{\phi( t-k)e^{\frac{-j}{2}[t^2-k^2 -(t-k)^2]\cot\alpha}}\,dt\right|^2\\\\
&=2^j\sum_{k \in \mathbb Z}\left|\dfrac{1}{2\pi\sin\alpha}\int_{-\pi\sin\alpha}^{\pi\sin\alpha} \hat{f}_\alpha(2^{j}u)\overline{\Theta_\alpha(u)}\,du\right|^2\\\\
&=\dfrac{2^j}{2\pi\sin\alpha}\int_{-\pi\sin\alpha}^{\pi\sin\alpha} \left|\hat{f}_\alpha(2^{j}u)\overline{\Theta_\alpha(u)}\right|^2\,du\\\\
&=\dfrac{1}{2\pi\sin\alpha}\int_{-2^j \pi\sin\alpha}^{2^j\pi\sin\alpha} \left|\hat{f}_\alpha(\xi)\overline{\Theta_\alpha(2^{-j}\xi)}\right|^2\,du\\\\
&\rightarrow\dfrac{1}{2\pi\sin\alpha}\int_{-\infty}^{\infty} \left|\hat{f}_\alpha(\xi)\right|^2\,d\xi =\|f\|_2^2
\end{align*}
as $j \rightarrow \infty$ because of the dominated convergence theorem.$\square$

\parindent=0mm \vspace{.1in}
\textbf{Remark:} The Characterization given in Theorem 5.1 shows that if $\phi$ is a scaling function for a fractional MRA, then the function $\Upsilon,$ defined by $\mathcal{F}_\alpha\{\Upsilon\}(u) = |\Theta_\alpha(u)|,$ is also a scaling function for fractional MRA of $L^2(\mathbb R)$.
\parindent=0mm \vspace{.1in}

{\bf{References}}

\end{document}